\documentclass{article}
\usepackage{multicol}
\usepackage{microtype}
\usepackage{theorem}
\usepackage{url}
\usepackage{algorithmicx}
\usepackage{algpseudocode}
\usepackage{enumerate}
\usepackage{amsfonts}
\usepackage{amsmath}
\usepackage{geometry}
\usepackage[english]{babel}
\newcommand{\FF}{\mathbb F}
\newcommand{\G}{\mathcal G}

\newcommand{\PG}{\mathrm PG}
\newcommand{\cP}{\mathcal P}

\newcommand{\cC}{\mathcal C}
\newcommand{\cD}{\mathcal D}

\newcommand{\codim}{\mathrm{codim}\,}

\newtheorem{mth}{Main Result}

\title{Polar Grassmannians and Their Codes}
\author{Ilaria Cardinali and Luca Giuzzi}
\date{December 12, 2014}
\begin{document}
\maketitle
\begin{abstract}
  We present a concise description of Orthogonal Polar Grassmann Codes
  and motivate their relevance. We also describe efficient encoding
  and decoding algorithms for the case of Line Grassmannians and
  introduce some open problems.
\end{abstract}
\section{Introduction}
The aim of this talk is to survey some recent results on Polar 
Grassmann codes and propose some open problems which we are currently
considering.

Polar Grassmann codes have been introduced in \cite{IL13} as the
projective codes arising from the Pl\"ucker embedding $\varepsilon_k$ of an
orthogonal Grassmannian $\Delta_{n,k}$.

The abstract is organized as follows: in Section \ref{emb} we provide some
quick background on orthogonal Grassmannians and describe
some of the projective codes
they determine.
Section \ref{ec} introduces an enumerator which might be
use to implement efficient encoding for line polar Grassmann codes.
Finally, in Section \ref{op} we enumerate some open problems which might
be of interest.

\section{Projective codes from polar Grassmannians}
\label{emb}

Let $V := V(2n+1,q)$ be a $(2n+1)$--dimensional vector space over a finite
field $\FF_q$ endowed with a non--singular quadratic form  $\eta$ of
Witt index $n$.

The \emph{polar Grassmannian of orthogonal type} $\Delta_{n,k}$ is a
proper subgeometry of the Grassmannian $\G_{2n+1,k}$ of
the $k$--subspaces of $V$ whose points correspond to all totally
singular $k$--spaces of $V$ with respect to $\eta$.

For $k<n$ the lines of $\Delta_{n,k}$ are exactly the
lines $\ell_{X,Y}$ of $\G_{2n+1,k}$ where $Y$ is totally singular; on the other hand,
when $k=n$ the
lines of $\Delta_{n}:=\Delta_{n,n}$ are sets of the form
\[\ell_{X}:=\{Z\mid X \subset Z\subset X^{\perp} , ~\dim(Z) = n,
~ Z ~\text{totally ~singular}\}\]
with $X$  a totally singular
$(n-1)$--subspace of $V$ and $X^\perp$ its orthogonal  with
respect to $\eta$. Note that the points of $\ell_{X}$ form a conic in the
projective plane $\PG(X^\perp/X)$.

Clearly, for $k=1$, the geometry $\Delta_{n,1}$ can be identified with
the orthogonal
polar space $Q(2n,q)$ of rank $n$ associated to $\eta$.
When $k=n$,
the geometry $\Delta_n$ is, in some sense, the dual of $\Delta_{n,1}$ and it
is called the
the \emph{orthogonal dual polar space} of rank $n$.

Define now $W_k :=\bigwedge^k V$. The \emph{Grassmann} 
 or \emph{Pl\"ucker}
 \emph{embedding}
 $e_k^{gr}:\G_{2n+1,k}\rightarrow \PG(W_k)$ 
 maps the arbitrary $k$--subspace
 $\langle v_1,v_2,\ldots, v_k\rangle$ of $V$ (hence a point of $\G_{2n+1,k}$)
 to the point $\langle
 v_1\wedge v_2\wedge\cdots \wedge v_k \rangle$ of $\PG(W_k)$.

It is well known that $e_k^{gr}(\G_{2n+1,k})$ is an algebraic variety;
furthermore the lines of $\G_{2n+1,k}$, that is to say collections of
$k$ spaces all mutually intersecting in a given $k-1$ space are mapped
into lines of $\PG(W_k)$.

Denote now by
$\varepsilon_k^{gr}:= {e_k^{gr}}|_{\Delta_{n,k}}$ be the truncated restriction of
$e_k^{gr}$ to $\Delta_{n,k}$.
The image of $\Delta_{n,k}$ under $\varepsilon_k$ is the $k$-Fano variety
of the quadric $Q(2n,\FF_q)$.
For $k<n$, all lines of $\Delta_{n,k}$ are mapped by $\varepsilon_k^{gr}$ into
lines of $\PG(W_k)$.  We observe however that for $k=n$ a line of 
$\Delta_n$ is mapped into a conic in $\PG(W_n)$; see \cite{CP1} for details.

Given any set of points $\Omega\subseteq\PG(W)$ with $W$ vector space
over $\FF_q$, the projective code defined
by $\Omega$ is a linear $q$--ary
code $\cC:=\cC(\Omega)$ whose generator matrix contains as columns
the coordinates of the points of $\Omega$. This code is uniquely defined up
to code equivalence. It is well known that the parameters of a projective
code are
\[ N=\#\Omega; \qquad K=\dim\langle\Omega\rangle;\qquad
 d_{\min}=\#\Omega-\max_{\footnotesize\begin{subarray}{c}\Pi\leq W,
    \Omega\not\subseteq\Pi\\
    \codim\Pi=1 \end{subarray}} \#(\Pi\cap\Omega).
\]
In particular, the study of the minimum distance of $\cC$ is equivalent
to the investigation of the possible hyperplane sections of $\Omega$.

Projective codes arising from the Pl\"ucker embedding $e_k$ of ordinary
Grassmannians
$\G_{n,k}$ in $\bigwedge^k V$
have been extensively investigated in recent years;
in particular, see \cite{R1,R2,N96,GL2001,GPP2009,GK2013}.

In the present note we are concerned with codes $\cP_{n,k}$ arising from the
projective system $\Omega_{n,k}:=\varepsilon_k(\Delta_{n,k})$ given by
the Pl\"ucker embedding $\varepsilon_k$ of an orthogonal Grassmannian.

\begin{mth}[\cite{IL13}]\label{main1}
  Let $\cP_{k,n}$ be the code arising from the projective
  system $\varepsilon_k^{gr}(\Delta_{n,k})$ for $1\leq k< n$.
  Then, the parameters of $\cP_{n,k}$ are
      \[N=
      \prod_{i=0}^{k-1}\frac{q^{2(n-i)}-1}{q^{i+1}-1},
      \qquad K=\left\{\begin{array}{ll}
          \binom{2n+1}{k} & \mbox{for $q$ odd} \\
          \binom{2n+1}{k}-\binom{2n+1}{k-2} &
          \mbox{for $q$ even,}\\
          \end{array}\right. \,\,\,\, \]
        \[ d\geq \psi_{n-k}(q)(q^{k(n-k)}-1)+1, \]
where $\psi_{n-k}(q)$ is the maximum size of a (partial) spread of the
parabolic quadric $Q(2(n-k),q)$.
\end{mth}
Here $\psi_r(q)=q^{r+1}+1$ for $q$ even and $\psi_r(q)\geq q+1$ for
$q$ odd.

To provide a sketch of the proof, observe that $N$ is just the number
of totally singular $k$--spaces contained in a $Q(2n,q)$.
The dimension $K$ arises from some results on embedding of Polar
Grassmannians in \cite{CP2,CP1}. In particular, for $q$ odd the dimension
of a polar Grassmann code is the same as that of the corresponding 
Grassmann code.
The estimate on the minimum distance derives from the study of maximal
totally singular subspaces contained in $Q(2n,q)$; in particular, for any
given fixed $k$--dimensional totally singular subspace $E$ there are at
least $\psi_{n-k}(q)$ generators $H_i$ of $Q(2n,q)$ meeting just in $E$.
In each of the spaces $H_i/E$ it is then possible to apply the well--known
result on the minimum distance of Grassmann codes.

As the sketch above illustrates, it has to be expected that the buonds
we obtain in Main Result \ref{main1} are not sharp.

More recently, in \cite{ILP14}, together with A. Pasini, we have
been able to fully determine the minimum distance for Line Polar Grassmann
codes, i.e. codes with $k=2$, for $q$ odd.
\begin{mth}[\cite{ILP14}]
\label{main3}
Suppose $n\geq 2$ and $q$ odd; then,
the minimum distance $d_{\min}$ of the orthogonal Grassmann code $\cP_{n,2}$
  is
  \[ d_{\min}=q^{4n-5}-q^{3n-4}. \]
  Furthermore, all words of minimum weight are projectively equivalent.
\end{mth}
The proof of this theorem hinges upon the observation that any hyperplane
$\Pi$ of $W_2=V\wedge V$ corresponds to an alternating bilinear form $\pi$
on $V$. 
In particular, we have $L\in\Pi\cap\varepsilon_2^{gr}(\Delta_{n,2})$ if,
and only if, the line $\ell$ whose image is $L$ is simultaneously totally
singular for the quadratic form $\eta$ and totally isotropic for the
(degenerate) bilinear form $\pi$. As it might be expected, the forms $\pi$
giving maximum intersection turn out to have maximum radical and,
consequently, maximum number of totally isotropic lines; in the case of
ordinary Grassmann codes, all of these forms are equivalent; however, for
polar Grassmann codes there are many inequivalent possibilities. Our 
argument relies upon providing a detailed analysis of these possibilities
and bounds on the values which might occur.

Observe that in Main Theorem \ref{main1} we have not considered the case of
dual polar spaces. We have however determined the value of the minimum distance
when $n=k=2$ and $n=k=3$, as illustrated by the following theorem.
\begin{mth}[\cite{IL13}]
\begin{enumerate}[(i)]
    \item\label{mt2:i}
      The code $\cP_{2,2}$ arising from a dual polar space of rank $2$ has parameters
      \[N=(q^2+1)(q+1),\qquad
      K=\left\{\begin{array}{ll}
          10 & \mbox{for $q$ odd} \\
          9  & \mbox{for $q$ even,}
          \end{array}\right. \qquad d=q^2(q-1).
       \]
    \item\label{mt2:ii}
      The code $\cP_{3,3}$ arising from a dual polar space of rank
      $3$ has parameters
      \begin{small}
      \[
      \begin{array}{ll l}
        N=(q^3+1)(q^2+1)(q+1),\,\,\,& K=35,\,\,\, & d=q^2(q-1)(q^3-1)\,\,\,  \mbox{ for $q$ odd} \\
         & and& \\
        N=(q^3+1)(q^2+1)(q+1),\,\,\,& K=28,\,\,\, & d=q^5(q-1)\,\,\, \mbox{ for $q$ even}.\\
       \end{array}\]
       \end{small}
    \end{enumerate}
\end{mth}

\section{Enumerative encoding}
\label{ec}
Grassmann linear codes have a very low data rate; as such it is paramount to be
able to describe efficient encoding and decoding algorithms acting
locally on the components.
To this aim, in \cite{SE}, an efficient algorithm for 
 enumerative coding of Grassmannians is introduced; see also \cite{M12}
for some improvement. Their work is based upon the approach of \cite{Cover}
which requires to determine the number of subspaces whose representation begins
with a given prefix.

The techniques of \cite{SE} cannot be directly applied to polar Grassmannians,
as the value of the quadratic form $\eta$ must also be tracked.
In \cite{IL14} we introduced an enumerator algorithm for Line Polar Orthogonal 
Grassmannians of complexity $O(q^2n^3)$.

Using this algorithm, we can provide both efficient encoding and efficient
error correction for Polar Grassmann codes. More in detail,
\begin{enumerate}
  \item It is possible to fully locally encode any Line (Polar) Grassmann code.
    Indeed, given a message $\mathbf{m}$ it is easy to determine an alternating
    form $m(x,y)$ acting on the vector space $V$. For any position $i$ in
    $\cP_{n,2}$, the value of the codeword corresponding to $\mathbf{m}$
    is just $m(A,B)$ where $A$, $B$ are the two generators of the line
    with index $i$ in $Q(2n,q)$ taken in Row Reduced Echelon Form.
  \item There is also a form of local error correction which can be
    obtained by exploiting the geometry. Indeed, given an index position
    $i$, let $\ell_i$ be the corresponding totally singular line.
    Then, it is possible to study the bilinear forms induced by a
    codeword $\mathbf{c}$ on the totally singular planes passing through
    $\ell_i$ and use this information in order to recover the value $c_i$
    is supposed to have.
\end{enumerate}
The details are contained in \cite{IL14}.

\section{Open problems}
\label{op}
We conclude this survey, by presenting some open problems.
\begin{enumerate}[a)]
\item \emph{Bounds for the minimum distance of Orthogonal
    Polar Gra{ss}mann codes for $k>2$.} \\
\item \emph{Generating sets of minimum weight} \\
  In the case of line polar Gra{ss}mann codes the minimum weight codewords
  are all projectively equivalent. Do they constitute a
  generating set for the code ? If not, what is the dimension of the
  subcode they span? 
\item \emph{Spectrum of low weight codewords} \\
  In \cite{ILP14} we constructed several classes of codewords (depending on
  some parameters) in order to explicitly estimate the minimum distance.
  The same construction can be used to produce several codewords of
  different weight.  Do they provide an exhaustive list of all possible
  small weight codewords? What about the weight enumerator?
\item \emph{Higher weights} \\
  Determine at least some of the higher weights of
  line polar Gra{ss}mann codes.
\item \emph{LDPC codes} \\
  Consider the incidence matrix $H$ of the design of
  hyperplanes of $\bigwedge^k V$ and
  the Gra{ss}mann embedding of a polar Gra{ss}mannian
  $\cD:=\varepsilon_k(\Delta_{n,k})$.
  What can we say about the binary code $\cC$ with parity check/generator
  matrix $H$?
  Some low weight codewords of $\cC$ correspond to hyperplanes with
  \emph{maximum} intersection with $\cD$; are these
  the codewords of minimum weight? 
  Is this code generated by its minimum  weight codewords?
  What is the dimension of $\cC$ (i.e. the $2$--rank of $G$) ?
  What about $p$--ranks for $p>2$?
\item \emph{Network coding with polar Grassmannians} \\
  The Grassmann graph arising from Polar Grassmannians
  has diameter strictly larger than that of the corresponding projective 
  Grassmannian. Can this be exploited in order to offer better correction
  capabilities in the case of random network coding?
\end{enumerate}

\end{document}